\documentclass{amsart}

\usepackage{amssymb}
\usepackage{stmaryrd}
\usepackage{enumerate}
\usepackage{hyperref}
\usepackage{mathrsfs}
\usepackage[matrix,arrow,graph,curve]{xy}

\theoremstyle{plain}      \newtheorem*{theorem}{Theorem}
\theoremstyle{plain}      \newtheorem*{lemma}{Lemma}
\theoremstyle{plain}      \newtheorem*{proposition}{Proposition}
\theoremstyle{plain}      
\theoremstyle{definition} 
\theoremstyle{definition} 
\theoremstyle{remark}     
\theoremstyle{remark}     \newtheorem*{note}{Note}

\renewcommand{\mathcal}{\mathscr}
\newcommand{\A}{\ensuremath{\mathcal{A}}}
\newcommand{\B}{\ensuremath{\mathcal{B}}}
\newcommand{\C}{\ensuremath{\mathcal{C}}}
\newcommand{\I}{\ensuremath{\mathcal{I}}}
\newcommand{\M}{\ensuremath{\mathcal{M}}}
\renewcommand{\P}{\ensuremath{\mathcal{P}}}
\newcommand{\X}{\ensuremath{\mathcal{X}}}
\newcommand{\V}{\ensuremath{\mathcal{V}}}
\newcommand{\Set}{\ensuremath{\mathbf{Set}}}
\newcommand{\GSet}{\ensuremath{\mathrm{G}\text{-}\Set}}
\newcommand{\Joy}{\ensuremath{\mathbf{Joy}}}
\newcommand{\Rep}{\ensuremath{\mathbf{Rep}}}
\newcommand{\Span}{\ensuremath{\mathbf{Span}}}
\newcommand{\Vect}{\ensuremath{\mathbf{Vect}}}

\newcommand{\x}{\ensuremath{\times}}
\renewcommand{\o}{\ensuremath{\circ}}
\newcommand{\ox}{\ensuremath{\otimes}}
\newcommand{\bigox}{\ensuremath{\bigotimes}}
\newcommand{\op}{\ensuremath{\mathrm{op}}}
\newcommand{\ob}{\ensuremath{\mathrm{ob}}}
\newcommand{\ol}{\ensuremath{\overline}}
\newcommand{\ul}{\ensuremath{\underline}}
\newcommand{\olstar}{\ensuremath{\mathbin{\overline{*}}}}
\newcommand{\ulstar}{\ensuremath{\mathbin{\underline{*}}}}
\newcommand{\oK}{\ensuremath{\overline{K}}}
\newcommand{\uK}{\ensuremath{\underline{K}}}
\newcommand{\<}{\ensuremath{\langle}}
\renewcommand{\>}{\ensuremath{\rangle}}
\newcommand{\dint}{\ensuremath{\displaystyle{\int}}}

\newcommand{\ra}[1][{}]{\ensuremath{\xymatrix@=15pt@1{\ar[r]^{#1}&}}}
\newcommand{\Ra}{\ensuremath{\xymatrix@C=15pt@1{\ar@{=>}[r]&}}}

\title{Monoidal functor categories and graphic Fourier transforms}
\author{Brian J. Day}
\address{Centre of Australian Category Theory, Macquarie University, NSW,
2109, Australia}
\date{15 December 2006}

\begin{document}

\begin{abstract}
This article represents a preliminary attempt to link Kan extensions, and
some of their further developments, to Fourier theory and quantum algebra
through $*$-autonomous monoidal categories and related structures.
\end{abstract}

\maketitle

\section*{Introduction}

In Part~1 of this article we use $*$-autonomous symmetric monoidal
categories $\V$ (in the sense of~\cite{1}) and their resulting functor
categories $\P(\A) = [\A,\V]$, to describe aspects of the ``graphic''
upper and lower convolutions of functors into $\V$, and their transforms.

Thus the obvious notion of ``Fourier'' transform (via Kan extension) of a
functor in a monoidal functor category $[\A,\V]$ is reiterated in \S 1.3
below, where the transform of the convolution product of two such functors
is seen to be the (often pointwise) tensor product of their transforms.

The basic notion of multiplicative kernel is defined in \S 1.2, this being
the main source of multiplicative functors. Examples are described, with a
brief account of the association scheme example at the end of \S 1.2.

In \S 1.3 we look at transforms of functors (called analytic functors after A.
Joyal~\cite{12}) in the context of what we call the Joyal-Wiener
category. This category of transforms is, under simple conditions, a
monoidal category equivalent to the original monoidal domain $[\A,\V]$; of
course, in Joyal's case, it is also equivalent to the category of all
Joyal-analytic functors and weakly cartesian maps between them.

In general, the transform of a functor behaves like a classical Fourier
transform, and there is usually a corresponding inversion process. For
instance it is possible (see \S 2.1) that a transformation functor can have
a tentative left inverse which then leads to the construction of the
required (two-sided) inversion data.

These transforms also generalize the Fourier transforms of Hopf algebras (as
discussed in~\cite{4} for example) which extend directly from Hopf algebras
in $\V$ to the many-object Hopf algebroids of~\cite{10}; some mention of this
is made in \S 2.1. Several other types of examples are described in \S 2.2.

The original convolution construction on a functor category of the form
$[\A,\V]$ for a small promonoidal structure $\A$, may be found in~\cite{6}
and~\cite{7}, where such categories are viewed in much the same light as
function algebras. It is emphasised here that all the categorical concepts
used below, such as ``category'', ``functor'', ``natural transformation'',
etc., are $\V$-enriched (in the sense of~\cite{13}) over the given symmetric
monoidal closed base category $\V$, unless otherwise mentioned in the text.

\section{Multiplicative kernels}

\subsection{Upper and lower convolution}

Let $\V = (\V_0,I,\ox,[-,-],(-)^*)$ be a complete (hence cocomplete)
$*$-autonomous symmetric monoidal closed category (in the sense
of~\cite{1}) with $I$ as dualising object. Recall that
\[
[X,Y] \cong (X \ox Y^*)^*,
\]
and if an object $Z$ has a dual $Z^\vee$ in $\V$, then
\[
Z^\vee \cong [Z,I] = Z^*
\]
in which case
\[
[Z,X] \cong Z^* \ox X \quad \text{and} \quad (Z \ox X)^* \cong Z^* \ox X^*
\]
for all $X$ in $\V$.

If $(\A,p,j)$ is a small promonoidal category over $\V$, then the \emph{upper
convolution} of $f$ and $g$ in the functor category $[\A,\V]$ is defined
in~\cite{6} as
\[
f \olstar g = \int^{ab} f(a) \ox g(b) \ox p(a,b,-) \quad \text{and} \quad
    \ol{I} = j.
\]

If $(\A^\op,p,j)$ is a small promonoidal category, then the \emph{lower
convolution} of $f$ and $g$ in $[\A,\V]$ is defined as
\[
f \ulstar g = (\int^{ab} f(a)^* \ox g(b)^* \ox p(a,b,-))^*
    \quad \text{and} \quad \ul{I} = j^*.
\]
(See~\cite{8} and~\cite{16} for example.)

Both products yield associative and unital monoidal structures on $[\A,\V]$;
the upper product preserves $\V$-colimits in each variable, while the lower
product preserves $\V$-limits in each variable. The upper product $f \olstar
g$ in $[\A,\V]$ using $p$ on $\A$ gives $f \olstar g$ on
\[
[\A,\V]^\op = [\A^\op,\V^\op]
\]
which transforms under the equivalence $\V^\op \simeq \V$ to the lower
product $f^* \ulstar g^*$ in $[\A^\op,\V]$ using the same $p$ on $\A$, since
\begin{align*}
(f \olstar g)^* &= (\int^{ab} f(a) \ox g(b) \ox p(a,b,-))^* \\
                &\cong (\int^{ab} f(a)^{**} \ox g(b)^{**} \ox p(a,b,-))^* \\
                &= f^* \ulstar g^*.
\end{align*}

An \emph{antipode} $S$ on a (promonoidal) category $(\A,p,j)$ is a functor
\[
S:\A^\op \ra \A
\]
such that $S^\op \dashv S$ with $S^2 \cong 1$. Then the resulting data
$q(a,b,c) = p(Sa,Sb,Sc)$ and $k(c) = j(Sc)$ become part of an obvious
promonoidal structure on $\A^\op$.

If the set of data $(\A,p,j,S)$ is an \emph{$S$-autonomous} category, in the
sense that the cyclic condition
\[
p(a,b,Sc) \cong p(b,c,Sa)
\]
holds naturally in $a,b,c \in \A$, then the derived set $(\A^\op,q,k,S^\op)$
is $S^\op$-autonomous since
\[
p(Sa,Sb,S^2 c) \cong p(Sb,Sc,S^2 a);
\]
that is
\[
q(a,b,Sc) \cong q(b,c,Sa).
\]

\begin{theorem}
If $(\A,p,j,S)$ is a $S$-autonomous promonoidal category and
$(\V,I,\ox,(-)^*)$ is the $*$-autonomous base category, then the upper
convolution structure on $[\A,\V]$ is $*$-autonomous under the antipode
defined by $f^*(a)=f(Sa)^*$ for $f$ in $[\A,\V]$.
\end{theorem}

\begin{proof}
We have the natural isomorphisms
\begin{align*}
[f,g](c) &= \int_{ab} [f(a) \ox p(c,a,b),g(b)] \text{ by definition of
$[-,-]$ in $[\A,\V]$,} \\
    &\cong \int_{ab} (f(a) \ox p(c,a,b) \ox g(b)^*)^* \text{ since $\V$ is
$*$-autonomous monoidal,} \\
    &\cong (\int^{ab} f(a) \ox g(b)^* \ox p(c,a,b))^*, \\
    &\cong (\int^{ab} f(a) \ox g(Sb)^* \ox p(c,a,Sb))^*
        \text{ since $S^2 \cong 1$,}\\
    &\cong (\int^{ab} f(a) \ox g^*(b) \ox p(a,b,Sc))^* \\
& \text{\quad since $p(c,a,Sb)=p(a,b,Sc)$ because $(\A,p,j,S)$ is $S$-autonomous,} \\
    &= (f \olstar g^*)^* \text{by definition of $\olstar$ on $[\A,\V]$}.
\end{align*}
\end{proof}

The upper convolution $f \olstar_p g$ is then related to the lower
convolution $f \ulstar_q g$ on the functor category $[\A,\V]$ using
$q(a,b,c)=p(Sa,Sb,Sc)$ on $\A^\op$, the latter being naturally isomorphic to
the ($*$-autonomous) product
\[
(f^* \olstar_p g^*)^*
\]
on $[\A,\V]$. This follows from the convolution calculation
\begin{align*}
(f^* \olstar_p g^*)^*(c) &= (\int^{ab} f(Sa)^* \ox g(Sb)^* \ox p(a,b,Sc))^*
        \text{ by definition of $\olstar_p$}, \\
    &\cong\int_{ab} (f(a)^* \ox g(b)^* \ox p(Sa,Sb,Sc))^*
        \text{ using $S^2 \cong 1$,} \\
    &\cong (\int^{ab} f(a)^* \ox g(b)^* \ox q(a,b,c))^*
        \text{ by definition of $q$ in $\A^\op$,} \\
    &= f(a) \ulstar_q g \text{ by definition of $\ulstar_q$.}
\end{align*}

\begin{note}
In the sequel we shall not go into the particular $*$-autonomous aspects of
the theory in any great detail.
\end{note}

\subsection{Multiplicative kernels}

For the given base category $\V$, a \emph{multiplicative kernel} from a
promonoidal $\A$ to another promonoidal $\X$ is a $\V$-functor
\[
K:\A^\op \ox \X \ra \V
\]
together with two natural structure isomorphisms
\begin{align*}
\int^{yz} K(a,y) \ox K(b,z) \ox p(y,z,x) &\cong \int^c K(c,x) \ox p(a,b,c)
    \quad \text{and} \\
j(x) &\cong \int^c K(c,x) \ox j(c).
\end{align*}
A $\V$-natural transformation
\[
\sigma:H \Ra K:\A^\op \ox \X \ra \V
\]
is called \emph{multiplicative} if it commutes with the structure
isomorphisms of $H$ and $K$.

A simple calculation with coends shows that the ``module composite'' of two
multiplicative kernels is again a multiplicative kernel and, for the given
$\V$, this leads to a monoidal bicategory, in the sense of~\cite{10}, with
promonoidal $\V$-categories as the objects (0-cells), multiplicative kernels
$K:\A^\op \ox \X \ra \V$ as the 1-cells $\A \ra \X$, and the multiplicative
$\V$-natural transformations $\sigma$ as the 2-cells.

Here are a few routine examples:

\subsection*{Examples}

\begin{enumerate}[(a)]
\item In the case where $\A$ is monoidal and $\X$ is comonoidal, so that
\begin{align*}
p(a,b,c) &= \A(a \ox b,c) \\
j(c) &= \A(I,c)
\end{align*}
and
\begin{align*}
p(y,z,x) &= \X(y,x) \ox \X(z,x) \\
j(x) &= I
\end{align*}
the structure isomorphisms for a multiplicative $K:\A^\op \ox \X \ra \V$
reduce to isomorphisms
\begin{align*}
K(a,x) \ox K(b,x) &\ra[\cong] K(a \ox b,x) \\
I &\ra[\cong] K(I,x)
\end{align*}
by the Yoneda lemma.

\item If $\A$ and $\X$ are both monoidal then a $\V$-functor $\phi:\A \ra
\X$ is multiplicative (i.e.,
\begin{align*}
\phi(a) \ox \phi(b) &\ra[\cong] \phi(a \ox b) \\
I &\ra[\cong] \phi(I))
\end{align*}
if and only if the module
\[
K(a,x) = \X(\phi(a),x)
\]
is a multiplicative kernel, again by Yoneda.

\item For any $\A$ and $\X$ promonoidal and $\V$-functor $\psi:\X \ra \A$
the conditions for the module
\[
K(a,x) = \A(a,\psi(x))
\]
to be a multiplicative kernel are precisely the conditions for restriction
along $\psi$ to be a multiplicative functor:
\[
[\psi,1]:[\A,\V] \ra[] [\X,\V]
\]
when $[\A,\V]$ and $[\X,\V]$ are given their upper convolution tensor
products.

\item If $\A=\I$ (the identity $\V$-category), then $K:\X \ra \V$ is a
multiplicative kernel if and only if $K \olstar K \cong K$. In two trivial
cases, take $\V=(0 \le 1)$ (cartesian closed) and $\A=1$. First let $\X=M$
be a module over a ring $R$, and define
\[
p(x,y,z) =
\begin{cases}
1 \text{ iff } z=rx+sy \text{ for some } r,s \in R \\
0 \text{ else.}
\end{cases}
\]
Then $K:M \ra (0 \le 1)$ is a multiplicative kernel if and only if
$K^{-1}(1) \subset M$ is a $R$-submodule. Secondly, let $\X=X$ be a
convexity space with join $[x,y]$ of points $x,y \in \X$. Define
\[
p(x,y,z) =
\begin{cases}
1 \text{ iff } z \in [x,y] \\
0 \text{ else;}
\end{cases}
\]
then the multiplicative kernels correspond to the convex subsets of $X$. The
promonoidal structure on $\X$ has no identity in these two examples.

\item Association schemes~\cite{3}. Given a promonoidal category $(\A,p,j)$
that has a unit object $I$ to represent the identity $j$ (i.e., $j(a) \cong
\A(I,a)$) and given $\X$ of the form $\B^\op \ox \B$ with the usual
promultiplication corresponding to $\B$-bimodule composition ``$\circ$'',
any functionally $\A$-indexed family of functors
\[
M_a:\B^\op \ox \B \ra \V
\]
with $M_I=\hom \B$, yields the multiplicative kernel
\[
K:\A^\op \ox (\B^\op \ox \B) \ra \V
\]
(where $K(a,x,y) = M_a(x,y)$) if and only if there exists a natural
``structure'' isomorphism
\[
M_a \circ M_b \cong \int^c p(a,b,c) \ox M_c.
\]
Again this follows directly from the Yoneda lemma applied to the
multiplicative kernel criteria for this particular example.

By also defining an antipode
\[
T:\B^\op \ra \B
\]
on $\B$, one can incorporate the notion of transpose matrix into this
setting by defining
\[
M^T_a(x,y) = M_a(Ty,Tx).
\]

Finally, in the case of an association scheme, one has that (for $\V=\Set$)
the category $\B^\op \ox \B$ corresponds to the cartesian product $X \x X$
of a set $X$ with itself, while $\A$ is the discrete category with
objects the members of the given partition of $X \x X$, the cardinals of
the respective promultiplication values $p(a,b,c)$ being the structure
constants of the association scheme, and $j$ being represented by the
identity relation on $X$. Here $(\A,p,j)$ has $p(a,b,c) \cong
p(b^*,a^*,c^*)$ under the antipode $Sa=a^*$ (the reverse relation on $a$),
while $T$ is just the identity function on $X$.

\end{enumerate}

\subsection{Transforms and analytic functors}

Given $K$ as before, define the \emph{($K$-)transform} of $f$ in $[\A,\V]$
as the (left) Kan extension~\cite{13}
\[
\oK(f)(x) = \int^a K(a,x) \ox f(a),
\]
and similarly the \emph{dual} transform of $h$ in $[\A^\op,\V]$ is defined
as the (right) Kan extension
\[
\oK^\vee (h)(x) = \int_a [K(a,x),h(a)]
\]
--- especially used here when $\V$ is $*$-autonomous, in which case we have
\begin{align*}
\oK^\vee (h) &= \int_a (K(a,x) \ox h(a)^*)^* \\
    &\cong (\int^a K(a,x) \ox h(a)^*)^* \\
    &= \oK(h^*)^*.
\end{align*}

\begin{theorem}
If $K$ is multiplicative, then
\begin{enumerate}[{\upshape (i)}]
\item $\oK$ preserves upper convolution, and 
\item $\oK^\vee$ preserves lower convolution.
\end{enumerate}
\end{theorem}

\begin{proof} \quad
\begin{enumerate}[{\upshape (i)}]
\item
\begin{align*}
\oK(f \olstar g) &= \oK(\int^{ab} f(a) \ox g(b) \ox p(a,b,-)) \\
  &= \int^c K(c,-) \ox \int^{ab} f(a) \ox g(b) \ox p(a,b,c) \\
  &\cong \int^{ab} \int^{yz} K(a,y) \ox K(b,z) \ox p(y,z,-) \ox f(a) \ox g(b) \\
  &\cong \int^{ab} \int^{yz} K(a,y) \ox f(a) \ox K(b,z) \ox g(b) \ox p(y,z,-) \\
  &= \int^{yz} (\int^{a} K(a,y) \ox f(a)) \ox (\int^b K(b,z) \ox g(b)) \ox p(y,z,-) \\
  &= \oK(f) \olstar \oK(g).
\end{align*}

\item For functors $h=f^*$ and $k=g^*$ in $[\A^\op,\V]$, we have
\begin{align*}
\oK^\vee (h \ulstar k) &\cong \oK((h \ulstar k)^*)^* \\
    &\cong \oK(f \olstar g)^* & \text{since } h \ulstar k = f^* \ulstar g^*
\cong (f \olstar g)^*, \\
    &\cong (\oK(f) \olstar \oK(g))^* & \text{since $\oK$ preserves
$\olstar$ by (i),} \\
    &\cong \oK(f)^* \ulstar \oK(g)^* \\
    &\cong \oK(h^*)^* \ulstar \oK(k^*)^* \\
    &= \oK^\vee (h) \ulstar \oK^\vee (k).
\end{align*}

\end{enumerate}
\end{proof}

Each $\V$-functor $K:\A^\op \ox \X \ra \V$ yields the standard Kan
adjunction
\[
(\epsilon,\eta):\oK \dashv \uK:[\X,\V] \ra[] [\A,\V],
\]
where
\[
\uK(f)(a) = \int_x [K(a,x),f(x)].
\]
In the case where $K$ is a multiplicative kernel we call $\oK$ an
``analytic'' or ``Fourier'' transformation if and only if $\eta$ is an
equaliser (i.e., a regular monomorphism) in $[\A,\V]$.

\begin{note}
Many of the examples in Part 2 have $\oK$ fully faithful which means
$\eta$ is an isomorphism, i.e., 
\[
[f,g] \cong [\oK(f),\oK(g)].
\]
Then, if $\V$ is $*$-autonomous, there results a dual isomorphism for
\[
\<f,g\> = \int^a f(a)^* \ox g(a),
\]
namely
\[
\<f,g\> \cong \<\oK(f),\oK(g)\>.
\]
This follows from
\begin{align*}
\<f,g\>^* &= (\int^a f(a)^* \ox g(a))^* \\
    &= \int_a (g(a) \ox f(a)^*)^* \\
    &= \int_a [g(a),f(a)] \\
    &= \int_x [\oK(g)(x),\oK(f)(x)] & \text{since $\oK$ is fully
faithful,} \\
    &= \int_x (\oK(g)(x) \ox \oK(f)(x)^*)^* \\
    &= (\int^x \oK(f)(x)^* \ox \oK(g)(x))^* \\
    &= \<\oK(f),\oK(g)\>^*
\end{align*}
which is a type of Parseval relation.
\end{note}

Returning to the case where $\eta$ is merely a regular monomorphism, we have
that
\[
\xymatrix{f \ar[r]_-\eta & \uK\oK(f)
\ar@<3pt>[r]^-{\eta \uK\oK}
\ar@<-3pt>[r]_-{\uK\oK \eta} & \uK\oK\uK\oK(f)}
\]
is an equaliser diagram in $[\A,\V]$ (see~\cite{2} and~\cite{14} for example).

For each such kernel $K$ one can construct a ``Joyal-Wiener'' category, here
denoted $\Joy(K)$, as follows. A map
\[
\alpha:\oK(f) \Ra \oK(g)
\]
in $[\X,\V]$ is called \emph{regular} when
\[
\oK\uK(\alpha) \oK(\eta) = \oK(\eta)\alpha;
\]
then using the equalizer hypothesis on $\eta$, each such regular $\alpha$
equals $\oK(\beta)$ for a unique $\beta:f \Ra g$ in $[\A,\V]$. With this
in mind, the $\V$-category $\Joy(K)$ is defined to be that subcategory of
the Kleisli category $[\A,\V]_T$ for the monad
\[
T = (\uK\oK,\uK \epsilon \oK,\eta),
\]
with the same objects as $[\A,\V]_T \subset [\X,\V]$ but with the equaliser
equations
\[
\xymatrix@R=20pt{
\Joy(K)(f,g) \ar[r] & [\X,\V](\oK(f),\oK(g)) \ar[r] \ar[dr]_-{(1,\oK(\eta)}
    & [\X,\V](\oK\uK\oK(f),\oK\uK\oK(g)] \ar[d]^-{(\oK \eta,1)} \\
&& [\X,\V](\oK(f),\oK\uK\oK(g))}
\]
in $\V$ defining its respective ($\V$-enriched) homs.

As a result, we obtain the usual Kleisli factorisation
\[
\xymatrix@R=15pt@C=5pt{
[\A,\V] \ar@{^{(}->}[rr] && [\X,\V] \\
& [\A,\V] \ar[ul]^-{\oK_T} \ar[ur]_-\oK},
\]
where $\oK_T$ is conservative, and $\Joy(K)$ is a $\V$-category with an
equivalence
\[
[\A,\V] \ra[\simeq] \Joy(K)
\]
and a conservative embedding into $[\X,\V]$.

To relate this to the work of Joyal~\cite{12}, we now consider the special
case where $\A$ is promonoidal and $\X=\V$ itself (monoidal under $\ox$).
Denote a fixed kernel $E$ by
\[
X(a) = E(a,X)
\]
and think of the coend
\[
F(X) = \int^a f(a) \ox X(a)
\]
as being an ``E-analytic'' functor of $X \in \V$ with ``coefficients'' $f$
in $[\A,\V]$. Thus we obtain obvious types of extensions of A. Joyal's
original notion of analytic functor~\cite{12}. That is, if we take $\A$ to be
the free $\V$-category on the groupoid \{finite sets and bijections\} and
let $E(a,X) = \bigox^a X$, then an $E$-analytic functor $F:\V \ra \V$ takes
the form
\begin{align*}
F(X) &= \int^a f(a) \ox X(a) & (X(a) = E(a,X)) \\
     &\cong \sum_a f(a) \underset{\text{Aut}(a)}{\ox} (\bigox\nolimits^a X).
\end{align*}
Note that, for $\V = \Set$ and $\V = \Vect_k$, the units $\eta_f$ of the
adjuction $\ol{E} \dashv \ul{E}$ are equalisers in $[\A,\V]$ because the
diagram
\[
\xymatrix{
f(b) \ar[r]^-{\eta_{f,b}} &
    \displaystyle{\int_X} [E(b,X), \displaystyle{\int_a} E(a,X) \ox f(a)]
        \ar[d]^-t \\
\displaystyle{\int_a} \A(a,b) \ox f(a) \ar[u]^-Y \ar[r]_-{\int^a m \ox 1}
    & \displaystyle{\int_a} E(a,b) \ox f(a)}
\]
commutes, where $t$ is the composite map
\begin{multline*}
\int_x [E(b,x),\int^a E(a,x) \ox f(a)] \ra[\text{proj}]
[E(b,b), \int^a E(a,b) \ox f(a)] \\
\ra[{[\text{id},1]}] [I,\int^a E(a,b) \ox f(a)] \cong \int^a E(a,b) \ox f(a),
\end{multline*}
and $\int^a m \ox 1$ is a monomorphism in $\V$ since
\[
m:\A(a,b) \ra E(a,b)
\]
is a cartesian natural monomorphism ($\A$ being a groupoid). Hence
\[
[\A,\V] \ra[\simeq] \Joy(E)
\]
so that, for $\V = \Set$, $\Joy(E)$ is equivalent to the category of all
Joyal-analytic functors on $\Set$ and weakly cartesian maps between them, as
one might expect. 

Finally, even in the case of general $\V$, one can define the convolution
product of two $E$-analytic functors
\[
F(X) = \int^a f(a) \ox X(a) \quad \text{and} \quad
G(X) = \int^a g(a) \ox X(a)
\]
by
\begin{align*}
F * G(X) &= \int^{ab} f(a) \ox g(b) \ox \int^c p(a,b,c) \ox X(c) \\
    &\cong \int^c (\int^{ab} f(a) \ox g(b) \ox p(a,b,c)) \ox X(c).
\end{align*}
The Hadamard product~\cite{5} can be defined also as
\[
F \times G(X) = \int^{a} (f(a) \ox g(a)) \ox X(a)
\]
if $\A$ has a comonoidal structure as well as its initial promonoidal
structure. In fact, if $K:\A^\op \ox \A \ra \V$ is multiplicative with
respect to these two structures, then the $K$-transform of $F$, defined by
\[
\oK(F)(X) = \int^b (\int^a f(a) \ox K(a,b)) \ox X(b),
\]
has the property
\[
\oK(F*G)(X) \cong (\oK(F) \x \oK(G))(X).
\]

\section{Examples of transformation functors}

\subsection{Examples where $\oK$ is fully faithful}

In this section we provide explicit inverses to various examples of
$\V$-functors of the form
\[
\oK:[\A,\V] \ra[] [\X,\V].
\]
In these examples, $\X$ may be a large $\V$-category so that the functor
category $[\X,\V]$ does not exist in the $\V$-universe; however the
underlying ``category'' $[\X,\V]_0$ makes some sense.

The method is to firstly find a left inverse $\Gamma$ to $\oK$, then show
that $\Gamma_0$ is faithful when restricted to the full image of $\oK_0$ in
$[\X,\V]_0$. For convenience, we shall here suppose that the unit $I \in
\V_0$ generates $\V_0$, so that $\Gamma$ is then $\V$-faithful on the full
image of $\oK$; by virtue of this monomorphism, the full image of $\oK$ is
usually a genuine $\V$-category. Thus, since
\[
\xymatrix{
[f,g]  \ar[rd]_-\cong \ar[r]^-\oK & [\oK(f),\oK(g)] \ar[d]^\Gamma \\
& [\Gamma \oK(f),\Gamma \oK(g)]}
\]
commutes, one has that $\Gamma$ is an isomorphism, hence that $\oK$ is fully
faithful. In other words
\[
[\A,\V] \ra[\simeq] \Joy(K)
\]
where $\Joy(K) \subset [\X,\V]$ is a full embedding, not merely a
conservative one.

In this part we shall suppose that $\V$ is a complete and cocomplete
symmetric monoidal closed category. Of course the upper transformation $\oK$
has a corresponding lower $\oK^\vee$ if $\V$ happens to be $*$-autonomous
also, in which case
\[
\oK^\vee(h) \cong \oK(h^*)^*
\]
so that $\Gamma \oK \cong 1$ if and only if
\[
\Gamma^\vee \oK^\vee(h) = \Gamma(\oK^\vee(h)^*)^* \cong h
\]
for all $h$ in $[\A^\op,\V]$.

\subsection*{Example 1}

One of the simplest examples of a Fourier transform is that obtained from a
Hopf algebroid in the sense of~\cite{10}. First $\A$ is given the promonoidal
structure
\[
p(a,b,c) = \A(a,Sb) \ox \A(b,c) \quad \text{and} \quad j(a) = I,
\]
where $S$ denotes the antipode of the algebroid, so that convolution on
$[\A,\V]$ becomes
\begin{align*}
f * g(c) &= \int^{a,b} f(a) \ox g(b) \ox p(a,b,c) \\
    &\cong \int^b \int^a f(a) \ox \A(a,Sb) \ox g(b) \ox \A(b,c) \\
    &\cong \int^b f(Sb) \ox g(b) \ox \A(b,c)
\end{align*}
on applying the Yoneda lemma.

Secondly, let $\X$ be $\A$ also; then the pointwise tensor on $[\X,\V]$ is
the usual
\[
f \ox g(c) = f(c) \ox g(c).
\]
If $K:\A^\op \ox \X \ra \V$ is taken to be the hom functor of $\A$, we then
get an isomorphism
\[
\oK:[\A,\V] \cong [\X,\V] \quad \text{with} \quad \oK(f*g) \cong f \ox g
\]
because, for such a Hopf algebroid, the so-called ``Fourier'' natural
isomorphism
\[
\A(a,Sb) \ox \A(b,c) \ra[\cong] \A(a,c) \ox \A(b,c)
\]
is always available. In the particular case where $\A$ is the single Hopf
algebra $H$, the Fourier isomorphism is the composite
\[
\xymatrix{
H \ox H \ar[d]_{1 \ox \delta} \ar[r]^\cong & H \ox H \\
H \ox H \ox H \ar[r]_{1 \ox S \ox 1} & H \ox H \ox H \ar[u]_{\mu \ox 1}}
\]
with inverse $(\mu \ox 1)(1 \ox \delta)$; see~\cite[\S2.3]{4}.

\subsection*{Example 2}

Given any (small) promonoidal $\V$-category $(\A,p,j)$, let $\A$ be $\A$
itself,
\[
K = p:\A^\op \ox \A^\op \ox \A \ra \V
\]
and $\X = \A^\op \ox \A$. Then
\[
\oK:[\A,\V] \ra[] [\A^\op \ox \A,\V]
\]
is given by
\[
\oK(f) = \int^a f(a) \ox p(a,-,-).
\]

Since $f*g = \int^{xy} f(x) \ox g(y) \ox p(x,y,-)$, we obtain
\begin{align*}
\oK(f*g)(u,v) &= \int^{xyz} f(x) \ox g(y) \ox p(x,y,z) \ox p(z,u,v) \\
    &\cong \int^z \int^x f(x) \ox p(x,z,v) \ox \int^y g(y) \ox p(y,u,z) \\
    &= \int^z \oK(f)(z,v) \ox \oK(g)(u,z) \\
    &= (\oK(f) \o \oK(g))(u,v),
\end{align*}
so $K$ is multiplicative. A tentative $\Gamma$ for this $\oK$ is given by
\[
\Gamma(F)(b) = \int^a F(a,b) \ox j(a)
\]
since
\begin{align*}
\Gamma \oK(f)(b) &= \int^a \oK(f)(a,b) \ox j(a) \\
    &= \int^a \int^x f(x) \ox p(x,a,b) \ox j(a) \\
    &\cong \int^x f(x) \ox \A(x,b) \qquad \text{since $p*j = \A(-,-)$} \\
    &\cong f(b) \qquad \text{by Yoneda.} \\
\end{align*}

\begin{proposition}
If $\V = \Vect_k$ and $j(d)$ is finite dimensional for all $d$, then
$\Gamma$ is $\V$-faithful on all of $[\A^\op \ox \A,\V]$ if $j:\A \ra \V$ is
$\V$-faithful and monomorphisms split in $[\A,\V]$.
\end{proposition}

\begin{proof}
We require $F \ra[] [j,F*j]$ to be monic for all $F$ in $[\A^\op \ox \A,\V]$.
But
\begin{align*}
[j,F*j](d,c) &= [j(d),F*j(c)] \\
    &= [j(d),\int^x F(x,c) \ox j(x)] \\
    &\cong \int^x F(x,c) \ox [j(d),j(x)]
\end{align*}
and
\begin{align*}
F(d,c) &\cong \int^x F(x,c) \ox \A(d,x) \quad \text{(Yoneda)} \\
    &\xymatrix@C=15pt@1{\ar@{>->}[r]&} \int^x F(x,c) \ox [j(d),j(x)]
\end{align*}
if $j$ is $\V$-faithful.
\end{proof}

\subsection*{Example 3}

Let $(\A,p,j)$ be a small promonoidal $\V$-category. A functor
\[
\phi:\A^\op \ra \V
\]
is called multiplicative if there are natural isomorphisms
\begin{align*}
\int^c \phi(c) \ox p(a,b,c) &\cong \phi(a) \ox \phi(b) \\
\int^c \phi(c) \ox j(c) &\cong I
\end{align*}
(in the sense of~\cite{9}). A natural transformation between two such
functors is multiplicative if it commutes with these isomorphisms.

The (naturally) comonoidal $\V$-category
\[
\X = \M(\A^\op,\V)
\]
is defined to be the free $\V$-category in the ordinary category of all
multiplicative functors from $\A^\op$ to $\V$, and multiplicative natural
transformations between them. The kernel
\[
K:\A^\op \ox \M(\A^\op,\V) \ra \V
\]
is given by evaluation $K(a,\phi) = \phi(a)$ so that
\[
\oK(f)(\phi) = \int^a f(a) \ox \phi(a).
\]

If each representable functor $\A(-,a)$ is multiplicative, we can define a
functor
\[
\Gamma:[\M(\A^\op,\V),\V] \ra[] [\A,\V]
\]
by
\[
\Gamma(F)(d) = F(\A(-,d)),
\]
so that $\Gamma \oK \cong 1$ by the Yoneda lemma. But this $\Gamma$ is
faithful on the full subcategory of $[\M(\A^\op,\V),\V]$ consisting of those
functors $F$ which preserve the Yoneda colimit
\[
\phi \cong \int^a \phi(a) \ox \A(-,a).
\]
Moreover, each $\oK(f)$ is clearly an $F$ with this property, so that
$\oK$ is a full embedding.

By using the assumption that each $\phi$ is multiplicative the kernel $K$
can be seen to be multiplicative.

\subsection*{Example 4}

Let $\C$ be a monoidal $\V$-category and let $\A \subset \C$ be a small
Cauchy dense promonoidal $\V$-subcategory. Let $\X = \C \ox \V$ with the
product monoidal structure from $\C$ and $\V$, and let
\[
K:\A^\op \ox (\C \ox \V) \ra \V
\]
be given by
\[
K(a,c,x) = \C(a,c) \ox x.
\]
Then 
\[
\oK(f)(c,x) = \int^a f(a) \ox \C(a,c) \ox x \quad \text{and} \quad
\Gamma(F)(a) = F(a,I)
\]
since
\begin{align*}
\Gamma \oK(f)(b) &= \oK(f)(b,I) \\
    &= \int^a f(a) \ox \C(a,b) \ox I \\
    &\cong \int^a f(a) \ox \A(a,b) \\
    &\cong f(b).
\end{align*}
Further application of the Yoneda lemma, and the calculus of coends, shows
that $\oK$ is multiplicative because
\begin{align*}
\int^a K(a,z,x) \ox p(b,c,a) &\cong \int^a (\C(a,z) \ox x) \ox \C(b \ox c,a)
\quad \text{since $p(b,c,a) = \C(b \ox c,a)$ on $\A$,} \\
    &\cong \C(b \ox c,z) \ox x \quad \text{since $\A \subset \C$ is Cauchy
dense,} \\
    &\cong \int^{z',z'',x',x''} (\C(b,z') \ox x') \ox (\C(c,z'') \ox x'')
\ox (\C(z' \ox z'',z) \ox [x'\ox x'',x], \\
    &= \int^{z',z'',x',x''} K(b,z',x') \ox K(c,z'',x'') \ox
p((z',x'),(z'',x''),(z,x))
\end{align*}
as required. Moreover, $\oK$ lands in the full subcategory of $[\C \ox
\V,\V]$ consisting of those $F$ for which the canonical map
\[
x \ox F(a,I) \ra[\cong] F(a,x)
\]
is an isomorphism, on which $\Gamma$ is clearly faithful.

\subsection*{Example 5}

Suppose $\V$ is $*$-autonomous and $\A$ is a small monoidal $\V$-category
with a given natural isomorphism
\[
\A(a,b) \cong \A(b,a)^*.
\]

Let $\X = [\A,\V]^\op$ (monoidal under convolution) and define
\[
K:\A^\op \ox [\A,\V]^\op \ra \V
\]
by $K(a,g) = g(a)^*$. Then 
\[
\oK(f)(g) = \int^a f(a) \ox g(a)^* = \langle g,f \rangle
\]
so that $\Gamma(F)(a) = F(\A(a,-))$ since
\begin{align*}
\Gamma \oK(f)(b) &= \int^a f(a) \ox \A(b,a)^* \\
    &\cong \int^a f(a) \ox \A(a,b) \\
    &\cong f(b).
\end{align*}

To show that $\oK$ is multiplicative, we first prove:

\begin{lemma}
$g(a)^* \cong [\A,\V](\A(a,-),g)^* \cong [\A,\V](g,\A(a,-))$ for all $g$ in
$[\A,\V]$.
\end{lemma}

\begin{proof}
\begin{align*}
\int_b [g(b),\A(a,b)] &= \int_b (g(b) \ox \A(a,b)^*)^* \\
    &\cong (\int^b g(b) \ox \A(b,a))^* \\
    &\cong g(a)^* \quad \text{by Yoneda.}
\end{align*}
\end{proof}

\begin{proposition}
$\oK$ is multiplicative.
\end{proposition}

\begin{proof}
We require a natural isomorphism
\[
\int^{g,h} K(a,g) \ox K(b,h) \ox [\A,\V](k,g*h) \cong \int^c K(c,k) \ox
p(a,b,c).
\]
By the lemma, the left side is isomorphic to
\begin{align*}
\int^{g,h} [\A,\V](g,&\A(a,-)) \ox [\A,\V](h,\A(b,-)) \ox [\A,\V](k,g*h)  \\
&\cong [\A,\V](k,\A(a,-)*\A(b,-)) \quad \text{by Yoneda,} \\
&\cong [\A,\V](\A(a \ox b,-),k)^* \quad \text{by the Lemma,} \\
&\cong k(a \ox b)^* \\
&\cong \int^c K(c,k) \ox p(a,b,c),
\end{align*}
as required.
\end{proof}

Also $\Gamma$ is faithful on the full image of $\oK$ since
\begin{align*}
\oK(f)(g)^* &\cong \int_a [f(a),g(a)] \\
    &\cong \int_a [\oK(f)(\A(a,-)),g(a)].
\end{align*}

\subsection*{Example 6}

Let $\V = \Vect_k$ for a fixed field $k$, and let $\X$ be a finite
promonoidal $\V$-category. Take
\[
\A = [\X,\Vect_f]^\op \quad \text{(which is monoidal under convolution)}
\]
and let
\[
K:\A^\op \ox \X \ra \V \quad \text{be evaluation.}
\]
Then
\[
\oK(f)(x) = \int^a f(a) \ox a(x) \cong f(\X(x,-)),
\]
and we choose
\[
\Gamma(F)(a) = \int_x[a(x),F(x)].
\]
Then
\begin{align*}
\Gamma \oK(f)(b) &= \int_x [b(x),\int^a a(x) \ox f(a)] \\
    &\cong \int^a \int_x b(x)^* \ox a(x) \ox f(a)
        \quad \text{since $\int^a$ is left exact,} \\
    &\cong \int^a \A(b,a) \ox f(a) \\
    &\cong f(b) \quad \text{by Yoneda.}
\end{align*}
Thus $\Gamma$ is faithful on the full image of $\oK$ since
\begin{align*}
\Gamma \oK(f)(\X(x,-)) &\cong f(\X(x,-)) \\
    &\cong \oK f(x);
\end{align*}
moreover, $\oK$ is multiplicative by Yoneda, as required.

\subsection*{Example 7}

Let $\C$ be a small braided compact closed category, let $\A = \C^\op$ with
the monoidal structure induced from $\C$, and let $\X$ be a (finite) Cauchy
dense full subcategory of $\A$ (cf.~\cite{11} and~\cite{15}).

Then $\X$ is promonoidal with respect to the structure induced by $\A$, namely
\begin{align*}
p(x,y,z) &= \C(z,x \ox y) \\
j(z) &= \C(z,I);
\end{align*}
the associative and unital axioms follow from the Cauchy density of $\X
\subset \A$. Moreover, here the unit of the $\oK \dashv \uK$ adjuction,
where $K(a,x) = \C(x,a)$, is an isomorphism since
\[
\eta_f:f(a) \ra \int_x [K(a,x), \int^b K(b,x) \ox f(b)]
\]
gives
\[
\xymatrix{f(a) \ar[r] \ar[dr]
    & \displaystyle{\int_x}[\C(x,a),\displaystyle{\int^b} \C(x,b) \ox f(b)]
        \ar[d]^-\cong \\
    & \displaystyle{\int_x}[\C(x,a) \ox f(x)]}
\]
by the Yoneda lemma, which becomes the isomorphism
\[
f(a) \ra[\cong] \int_x [\C(x,a),f(x)]
\]
by the Cauchy density of $\X \subset \A$.

Also, $\oK$ is multiplicative because the kernel $K$ satisfies the
conditions
\begin{align*}
\int^{xy} K(a,x) \ox K(b,y) \ox p(x,y,z)
&= \int^{xy} \C(x,a) \ox \C(y,b) \ox \C(z,x \ox y) \\
&\cong \C(z,a \ox b) \quad \text{by Cauchy density of $\X \subset \A$,}
\end{align*}
while
\begin{align*}
\int^c K(c,z) \ox \C(c, a \ox b)
&= \int^c \C(z,c) \ox \C(c,a \ox b) \\
&\cong \C(z,a \ox b)
\end{align*}
by the Yoneda lemma applied to $c \in \C$, and
\[
j(x) \cong \int^c \C(x,c) \ox j(c) = \int^c K(c,x) \ox j(c)
\]
for the same reason.

\subsection{When $\oK$ is conservative}

Sometimes $\oK$ is only conservative; i.e., the unit of the $\oK \dashv \uK$
adjunction is merely an equaliser, not an isomorphism. Then, as pointed out
earlier, we obtain an equivalence
\[
[\A,\V] \simeq \Joy(K)
\]
with a conservative multiplicative embedding of $\Joy(K)$ into $[\X,\V]$. In
fact, examples are readily available if we assume that all strong
monomorphisms in $\V$ are equalisers, and that the canonical identity map
\[
\textrm{id}:I \ra \int_x [K(a,x),K(a,x)]
\]
is a strong monomorphism, preserved by the $\ox$ in $\V$, for all $a \in
\A$: In addition, we only need that each coprojection
\[
K(a,x) \ox f(a) \ra \int^b K(b,x) \ox f(b)
\]
is a strong monomorphism in $\V$ and each $K(a,x)$ has a dual in $\V$. This
follows from the consideration of the diagram
\[
\xymatrix{f(a) \ar[r]^-{\textrm{id} \ox 1} \ar[d]_{\eta_f}
    & \displaystyle{\int_x} [K(a,x),K(a,x)] \ox f(a) \ar[d]^\cong \\
\displaystyle{\int_x} [K(a,x),\displaystyle{\int^b} K(b,x) \ox f(b)]
    & \displaystyle{\int_x} K(a,x)^* \ox K(a,x) \ox f(a) \ar[l]}
\]
which commutes by observing that both legs are natural in $f \in [\A,\V]$,
and then $f$ can be replaced by a representable functor and the Yoneda lemma
employed.

\subsection*{Example 8}

The coprojections
\[
K(a,x) \ox f(a) \ra \int^b K(b,x) \ox f(b)
\]
are coretractions if
\[
\A(a,b) \cong \A(b,a)^*
\]
naturally in $a,b \in \A$, where we suppose that the hom-spaces of $\A$
have duals in $\V$. This result is fairly immediate from the Yoneda lemma.

\subsection*{Example 9 (cf. \S 2.1 Example 2)}

Let $\C$ be a cocomplete monoidal category for which each functor of the
form
\[
- \ox C:\C \ra \C
\]
preserves colimits, and suppose that there exists a full embedding
\[
N:\A^\op \subset \C
\]
where $\{Nx; x \in \A\}$ is a small dense set of projectives in $\C$.
Suppose also that $I = NI$ for some $I \in \A$. 

Define the functor $U:\C \ra[] [\A^\op \ox \A,\V]$
by
\[
U(C)(x,y) = \C(Ny,Nx \ox C).
\]
If $[\A^\op \ox \A,\V]$ is given the monoidal structure of bimodule
composition, then $U$ becomes a conservative multiplicative functor and
there is also an induced monoidal equivalence
\[
\ol{N}:\C \simeq [\A,\V]
\]
given by $\ol{N}(C)(x) = \C(Nx,C)$, where $\A$ has the promonoidal structure
for which
\[
p:\A^\op \ox \A^\op \ox \A \ra \V \quad \text{corresponds to} \quad 
UN:\A^\op \ra[] [\A^\op \ox \A,\V]
\]
and
\[
j = \A(I,-):\A \ra \V,
\]
and where $[\A,\V]$ has the resulting convolution monoidal structure.

\subsection*{Example 10 \textrm{($\V =$ commutative monoids in $\Set$;
see~\cite{17})}}

Let $G$ be a finite group, let $\A=\Span(\GSet_f)$, and let $\X=\Rep_f(G)$
(for a fixed field) viewed as $\V$-categories. Define the kernel
\[
K:\A^\op \ox \X \ra \V
\]
by $K(a,x) = \X(Fa,x)$ where
\[
F:\A \ra \X = \Rep_f(G)
\]
is the extension of the canonical functor
\[
F:\GSet_f \ra \X,
\]
defined by ``summing over the fibres'' (see~\cite[\S 10]{17} for example),
and define the functor
\[
U:\X \ra \Span(\GSet)
\]
by $U(x)$ is the large induced $G$-Set of the representation $x:G \ra
\Vect_f$ in $\X$.

Then we have a canonical natural transformation
\[
\tau:\X(Fb,Fa) \ra \Span(\GSet)(b,UFa)
\]
in $\V$, and a natural coretraction
\[
\beta:a \ra UF(a)
\]
in $\Span(\GSet)$. The fact that $\oK$ is conservative now follows from
commutativity of the diagram
\[
\xymatrix{
f(a) \ar[r]^-{\eta_f} \ar[dd]_{\cong}
    & \dint_{\hspace{-1ex}x}[\X(Fa,x),\dint^b \X(Fb,x) \ox f(b)] \\
    & \dint^b \X(Fb,Fa) \ox f(b) \ar[u]_{\cong} \ar[d]^-{\int^b \tau \ox 1} \\
\dint^b \A(b,a) \ox f(b) \ar[r]_-\rho
    & \dint^b \Span(\GSet)(b,UFa) \ox f(b)}
\]
where $\rho$ is the coretraction induced by $\beta$, and the isomorphisms
shown in the diagram are again by the Yoneda lemma.

The main result of~\cite{9} applies to show that $\oK$ is multiplicative.

\subsection*{Example 11 \textrm{(cf. \S 1.2)}}

As mentioned in \S 1.2(c), an arbitrary $\V$-functor
\[
\psi:\X \ra \A,
\]
with $\A$ and $\X$ small promonoidal $\V$-categories, induces a
multiplicative functor
\[
\oK = [\psi,1]:[\A,\V] \ra[] [\X,\V]
\]
between the respective convolutions if and only if the kernel functor on
$\A^\op \ox \X$ given by
\[
K(a,x) = \A(a,\psi(x))
\]
is multiplicative.

Of course, $\oK = [\psi,1]$ is conservative if $\psi$ is a surjection on
object sets.

\subsection*{Example 12}

In Example 2 of \S 2.1 we have that
\[
\oK:[\A,\V] \ra{} [\A^\op \ox \A,\V]
\]
is conservative if $(\A,p,j)$ is a closed category; i.e., if
\[
p(a,x,y) = \A(a,[x,y]) \quad \text{and} \quad j(y) = \A(I,y).
\]
In fact, here $\eta_f$ is a coretraction since
\begin{align*}
\eta_f:f(a) \ra & \int_{xy} [p(a,x,y), \int^b p(b,x,y) \ox f(b)] \\
     \cong & \int_{xy} [\A(a,[x,y]), f([x,y])] \quad \text{by Yoneda},
\end{align*}
where
\[
\xymatrix{
f(a) \ar[r]^-{\eta_f} \ar[d]_1 & \dint_{\hspace{-1ex}xy} [\A(a,[x,y]),f([x,y])]
    \ar[d]^{x=I,y=z} \\
f(a) \ar[r]^-\cong & \dint_{\hspace{-1ex}z} [\A(a,z),f(z)]}
\]
commutes.

\section*{Appendix}

The term ``graphic'' applied to coends, convolution products, transforms,
etc., in the above context is intended to relate especially to base
categories like $\V = \Set$ and $\V=\Vect_k$, where there is a definite
notion of finite graph. In such cases we call a given $\V$-subgraph $\X$ of
the $\V$-category $\A$ a finite $\V$-graph if $\ob(\X)$ is finite and each
$\V$-object $\X(x,y)$ of edges is finite (or finite dimensional).

Then, if the category $\A$ is the directed union of all its finite
$\V$-subgraphs $\X_\phi$ (or some convenient subset of these), there results
a canonical isomorphism
\[
\underset{\phi}{\mathrm{colim}} \int^x T_\phi(x,x) \ra[\cong] \int^a T(a,a),
\]
where $T_\phi$ denotes the restriction of each $\V$-functor
\[
T:\A^\op \ox \A \ra \V
\]
to the $\V$-graph $\X_\phi^\op \ox \X_\phi$, and where each finite ``coend''
$\int^x T_\phi(x,x)$ over $x \in \X_\phi$ is computed in the same way as the
usual coend of a functor over any (small) category. In particular, any finite
$\V$-limit which commutes with each finite $\int^x$, also commutes with
their filtered colimit $\int^a$ over $a \in \A$. 

However, we note that for many practical purposes the finite ``coend''
$\int^x$ can be replaced by the (usual) coend over the corresponding full
subcategory of $\A$ determined by $\ob(\X)$.



\end{document}